\documentclass[10pt, a4paper]{amsart}

\usepackage[left=3cm,right=2.5cm,top=3cm,bottom=3cm]{geometry}
\usepackage[english]{babel}
\usepackage[utf8]{inputenc}
\usepackage{amsmath}
\usepackage{amssymb} 
\usepackage{dsfont} 
\usepackage{graphicx}
\usepackage[colorinlistoftodos]{todonotes}
\usepackage[arrow, matrix, curve, pic]{xy} 
\usepackage{hyperref}
\usepackage[user,titleref]{zref}
\usepackage{enumitem}
\usepackage[T1]{fontenc}
\usepackage{fancyhdr}
\usepackage[utf8]{inputenc}

\theoremstyle{plain}    
\newtheorem{thm}{Theorem}[section]
\newtheorem{prop}[thm]{Proposition}

\newtheorem{lem}[thm]{Lemma}

\newtheorem{fact}[thm]{Fact}

\theoremstyle{remark}

\newcommand{\T}{\mathcal{T}}

\newcommand{\Om}{\Omega}
\newcommand{\sub}[1]{\noindent\textbf{#1.}}
\newcommand{\eproof}{\hfill $\Box$}

\title{The Zariski-Lipman conjecture for log canonical spaces}

\author{Stefan Heuver}

\begin{document}

\begin{abstract}
In this paper we give an elementary proof of the Zariski-Lipman conjecture for log canonical spaces.
\end{abstract}

\maketitle

\tableofcontents

\section{Introduction}
The Zariski-Lipman conjecture claims that an $n$-dimensional complex variety with locally free tangent sheaf of rank $n$ is smooth. Although not proven in general, this conjecture holds in special cases (see \cite{Druel13}). In 2013 Druel has proven this conjecture for log canonical spaces by using foliations and the Camacho-Sad formula. In 2014 Graf-Kov\'acs obtained the result of Druel by strengthening the extension theorem (see \cite[Theorem 1.5]{GKKP11}) for $1$-forms on log canonical pairs (see \cite{GK14}). For more information on log canonical spaces, see \cite{Reid87} or \cite{KM98}. The goal of this paper is to give a more elementary proof for log canonical surfaces in the style of \cite{GKK10} and to conclude the Zariski-Lipman conjecture for log canonical spaces by using the reduction technique of Druel (see \cite[Theorem 5.2]{Druel13}).
\begin{thm}[The Zariski-Lipman conjecture for log canonical spaces]
Let $X$ be a log canonical variety of dimension $n$ such that the tangent sheaf $\T_X$ is locally free of rank $n$. Then $X$ is smooth.
\end{thm}
To prove the $2$-dimensional case, we will use an argument of \cite{vSS85}. The idea is that under the given properties, a smooth $1$-form on the variety will extent to a smooth $1$-form on the resolution, which leads to a contraction, unless the variety has already been smooth. After this we will reduce the $n$-dimensional case to the surface case by using hyperplane sections.\\

\sub{Conventions and basics}
In this paper a variety is a integral, separated scheme of finite type over an algebraic closed field. Varieties are also assumed to be reduced and irreducible. In \cite{Lip65} Lipman has proven that his conjecture fails to be true if $X$ is a variety over a field with positive characteristic and that a variety with locally free tangent sheaf is necessarily normal. That is why we will work over the field of complex numbers and assume that $X$ is normal. The sheaf $\Om^1_X$ of Kähler differentials behaves badly near a singular point. Therefore it is more useful to work with the reflexive hull $\Om^{[1]}_X$, which is the double dual of $\Om^1_X$ (see \cite[(1.5),(1.7)]{Reid87}). For the definition of logarithmic differential forms and logarithmic vector fields we recommend \cite{Saito80} and \cite {EV92}.\\

\sub{Acknowledgement}
The author would like to thank Daniel Greb for interesting discussions on this topic and great help with the master's thesis this paper is based on. 

\section{The 2-dimensional case}
To prove the Zariski-Lipman conjecture for log canonical surfaces, we use the minimal resolution $\pi:Y\rightarrow X$ and show that a reflexive $1$-form on $X$ lifts to a regular $1$-form on $Y$ under $\pi$.  An example of \cite[(1.8), (1.9)]{Reid87} shows that this is wrong if $X$ is not log canonical. In general a reflexive $1$-form only lifts to a $1$-form that is regular outside the exceptional locus. We will have to show that it extends regular over the exceptional locus.
\begin{prop}[Logarithmic extension]
Let $X$ be a log canonical surface with locally free tangent sheaf $\T_X$ of rank 2. Let $\omega\in H^0(X,\,\Om^{[1]}_X)$ be a reflexive $1$-form, $\pi:Y\rightarrow X$ the minimal resolution and $E$ the largest reduced divisor included in $excep(\pi)$. Then $\omega$ extends to a logarithmic $1$-form $\widetilde{\omega}:=\pi^*\omega\in H^0\left(Y,\,\Om^{1}_Y(\mathrm{log} E)\right)$ on $Y$.
\end{prop}
\proof
(cf. \cite[Proposition 6.1]{GKK10}) Let $\omega\in H^0(X,\,\Om^{[1]}_X)$. Since $\T_X$ is locally free of rank 2 we can assume without loss of generality that $\mathcal O_X( \mathcal K_X)\cong \mathcal O_X$. Since $\Om^{[1]}_X$ is the dual of $\T_X$ there exists a unique vector field $\xi\in H^0(X,\, \T_X)$ corresponding to $\omega$ via the perfect pairing
$$
 \Om^{[1]}_X\times\Om^{[1]}_X\rightarrow \mathcal O_X(\mathcal K_X)\cong \mathcal O_X.
$$
The minimal resolution $\pi$ is functorial and we can lift the vector field $\xi$ to a vector field $\widetilde{\xi}\in H^0(Y,\, \mathcal T_Y(-\mathrm{log} E))$. As $X$ is assumed to be log canonical we have $\mathcal O_Y(\mathcal K_Y + E)\cong \mathcal O_Y(D)$ for some effective divisor $D$ on $Y$. Hence $\widetilde{\xi}$ corresponds to an element $\widetilde{\omega}\in H^0(Y,\,\Om^1_Y(\mathrm{log} E)\otimes \mathcal O_Y(-D))$ via the pairing
$$
 \Om^{1}_Y(\mathrm{log} E)\times\Om^{1}_Y(\mathrm{log} E)\rightarrow \mathcal O_Y(\mathcal K_Y + E)\cong \mathcal O_Y(D).
$$
This yields the extension of $\omega$.
\eproof\\

The following result is a consequence of the negative definiteness of the self-intersection form in $E$.
\begin{prop}
Let $X$ be a normal surface and $\pi$ and $E$ as in Proposition 2.1. Then the inclusion
\[ H^0\left(Y,\,\Om^{1}_Y\right)\hookrightarrow H^0\left(Y,\,\Om^{1}_Y(\mathrm{log} E)\right)\]
is an isomorphism.
\end{prop}
\proof
See \cite[Lemma 1.3.b]{Wah85}.\eproof
\begin{thm}[The Zariski-Lipman conjecture for log canonical surfaces]
Let $X$ be a log canonical surface such that the tangent sheaf $\T_X$ is locally free of rank $2$. Then $X$ is smooth.
\end{thm}
\proof
Let $\pi$, $Y$, and $E$ be as defined above. Using Proposition 2.1 and 2.2 we see that a reflexive $1$-form on $X$ lifts to a regular $1$-form on $Y$ under $\pi$. With $\pi_*(\T_Y(- \mathrm{log} E))\cong \T_X$ Theorem 2.3 is a consequence of a classical argument presented in \cite[(1.6)]{vSS85}.\eproof
\section{The n-dimensional case}

In this section $X$ is an $n$-dimensional log canonical variety with locally free tangent sheaf and $\pi:Y\rightarrow X$ a functorial resolution (see \cite[3.45]{Kol07}). We will use $n-2$ hyperplane sections $G_1,\dots, G_{n-2}\subset X$ to cut $X$ down to a surface $S:=X\cap G_1 \cap \dots \cap G_{n-2}$ and show that $S$ is a log canonical surface with locally free tangent sheaf and therefore already smooth. Using the fact that a singularity of $X$ necessarily is a singularity of $S$, we conclude that $X$ must have been smooth. For more information on hyperplane sections used in this paper see \cite[2.E.]{GKKP11}. With $H:=\pi^{-1}(G)$ we will denote the preimage of $G$ under $\pi$. Please note the following facts. 
\begin{fact}[see {\cite[2.E.]{GKKP11}}]
Let $X$, $G$, $H$ and $\pi$ be as defined above and $E$ be the largest reduced divisor contained in $excep(\pi)$ then $G$ is affine and normal, $H$ is smooth and $\pi_{\mid H}$ is a functorial resolution and the largest reduced divisor $E_{\mid H}=E\cap H$ contained in $excep(\pi_{\mid H})$ is a simple normal crossing (snc) divisor (this uses Bertini's Theorem).
\end{fact}
\noindent Define $T:=Y \cap H_1 \cap \dots \cap H_{n-2}$, then by induction we get
\begin{enumerate}
	\item $S$ is affine and normal and $T$ is smooth,
    \item $\pi_{\mid T}:T\rightarrow S$ is a functorial resolution and the largest reduced divisor $C:=E_{\mid T}$ contained in the exceptional locus $excep(\pi_{\mid T})$ of $\pi_{\mid T}$ is a snc divisor. 
\end{enumerate}
We will now show that $S$ has the right properties.\\

\sub{The behavior of the singularities under reduction} 
The following lemma gives us a connection between $Sing(X)$ and $Sing(G)$.
\begin{lem}
Let $X$ be a normal variety and $G\subset X$ an effective, ample Cartier divisor. If G is smooth, then $Sing(X)\cap G=\emptyset$ and the singularities of $X$ are isolated. 
\end{lem}
\proof
This is a consequence of \cite[Lemma 1]{Che96} and the Jacobian criteria.\eproof
\begin{lem}
The surface $S$ constructed above is log canonical.
\end{lem}
\proof
Since $X$ is log canonical and normal, $|\mathcal G|$ is a basepoint free system of Cartier divisors and $G_1\in|\mathcal G|$ is a hyperplane section, Theorem 1.13 in \cite{Reid80} shows that $G$ is log canonical, too. Using Fact 3.1 we get Lemma 3.3 by induction.\eproof\\

\sub{The tangent sheaf of S is locally free}
To show that $\T_S$ is locally free we first need to prove the following Lemma.
\begin{lem}
Let $Y$ be a smooth variety of dimension $n\ge 2$ and $E\subset Y$ a snc divisor. Let $H\subset Y$ be a smooth hyperplane such that $E_{\mid H}$ is a snc divisor. Then the sequence
$$
0\rightarrow \mathcal N^*_{H\mid X} \rightarrow \Om^1_{Y}(\mathrm{log} E)_{\mid H} \rightarrow \Om^1_{H}(\mathrm{log} E_{\mid H}) \rightarrow 0
$$
is exact.
\end{lem}
\proof
Since $E$ and $E_{\mid H}$ are snc divisors we can use the sequence of \cite[2.3 a]{EV92}. Using the Snake Lemma we then get Lemma 3.4.\eproof
\begin{thm}
Let $X$ be a log canonical variety of dimension $n$ with locally free tangent sheaf of rank $n$. Then the tangent sheaf $\T_S$ of the surface $S$ defined above is locally free of rank $2$.
\end{thm}
\proof
The following proof is basically the cutting-down procedure of Druel (see the proof of \cite[Theorem 5.2]{Druel13}) supplemented with additional steps for the convenience of the reader. For $n=2$ the theorem is clear and we can assume that $n\ge 3$. Suppose that $Sing(X)\neq \emptyset$. Since $X$ is normal $codim_X(Sing(X))\ge 2$ and with \cite[p.318]{Flenner88} we get that $codim_X(Sing(X)) = 2$. Replacing $X$ with an affine open dense subset we may assume that $X$ is affine, $Sing(X)$ is irreducible of codimension $2$ and $\T_X\cong\mathcal O_X^{\oplus n}$.

Let $\pi: Y\rightarrow X$ be a functorial resolution and $E$ the largest reduced divisor contained in $excep(\pi)$. Note that $E\neq \emptyset$. We consider the morphism of vector bundles 
$$
F: \pi^*\T_X\rightarrow \T_Y(-\mathrm{log} E).
$$
Since $\T_X\cong\pi_*\T_Y(-\mathrm{log} E)$ the morphism $F$ is induced by the evaluation map
$$
\pi^*(\pi_*\T_Y(-\mathrm{log} E))\hookrightarrow \T_Y(-\mathrm{log} E)
$$
and induces an injective map of sheaves
$$
 \pi^*(\mathcal O_X(-\mathcal K_X))\cong \pi^* det(\mathcal T_X)\overset{\widetilde{F}}{\hookrightarrow} det(\mathcal T_{Y}(-\mathrm{log} E))\cong \mathcal O_{Y}(-\mathcal K_{Y}-E).
$$
Using the ramification formula $K_Y:=\pi^*K_X+\sum a_iE_i$ this yields $a_i\le -1$. Since $X$ is log canonical we get $a_i=-1$. Thus $F$ is an isomorphism. Since $\T_X$ is free we deduce that $\T_Y(-\mathrm{log} E)$ is free and that $\Om^1_{Y}(\mathrm{log} E)_{\mid H_1}\cong \mathcal O_{H_1}^{\oplus n}$.

Let $G_1\subset X$ be a general hyperplane section and $H_1=\pi^{-1}(G_1)\subset Y$. By Lemma 3.4 we have the exact sequence
$$
0\rightarrow \mathcal N^*_{H_1\mid X} \overset{\varPhi}{\rightarrow} \Om^1_{Y}(\mathrm{log} E)_{\mid H_1} \overset{\Psi}{\rightarrow}\Om^1_{H_1}(\mathrm{log} E_{\mid H_1}) \rightarrow 0.
$$
We want to prove
$$
\Omega^1_{H_1}(\mathrm{log} E_{\mid H_1})\cong \mathcal O_{H_1}^{\oplus dim(H_1)}.
$$

Since $X$ is affine and $G_1\in |\mathcal G|$ is an effective Cartier divisor and $|\mathcal G|$ basepoint free, we can define $G_1$ by a global function $g$. Due to this the ideal sheaf $\mathcal O_X(-G_1)=\mathcal N^*_{G_1\mid X}$ is free. Since $H_1$ is the total transform of $G_1$ we get $\mathcal N^*_{H_1\mid X}\cong\pi^*\mathcal N^*_{G_1\mid X}\cong \mathcal O_{H_1}$. Thus we can represent the map
$$
\varPhi: \mathcal O_{H_1}\cong\mathcal N^*_{H_1\mid X}\rightarrow \Om^1_{Y}(log E)_{\mid H_1}\cong \mathcal O_{H_1}^{\oplus n}
$$ 
by regular functions $f_1,\dots,f_n$ on $H_1$ and since $\pi_*(\mathcal O_{H_1})\cong \mathcal O_{G_1}$ by regular functions $g_1,\dots,g_n$ in $G_1$ with $f_i=g_i\circ\pi_{\mid H_1}$. If $g_{i\mid G_1\cap Sing(X)}=0$ for all $i$ then $\varPhi$ would vanish at every point of $\pi^{-1}(G_1\cap Sing(X))$. Since $E_{\mid H_1}$ is snc divisor, the sheaf $\Om^1_{Y}(\mathrm{log} E)_{\mid H}$ is locally free, which yields in a contradiction. Take an $i\in\{1,\dots, n\}$, so that $g_{i\mid G_1\cap Sing(X)}\neq 0$ then, by replacing $X$  with $X\backslash\{g_i=0\}$, we can ensure that $\pi_{\mid H_1}(x)\neq 0$ for all $x\in H_1$ and thus assume that $\varPhi$ has full rank $1$ in every fiber. We obtain the following exact sequence:
$$
0\rightarrow \mathcal O_{H_1}\overset{\varPhi}{\rightarrow} \mathcal O_{H_1}^{\oplus n}\overset{\Psi}{\rightarrow}\mathcal O^{\oplus(n-1)}_{H_1}\rightarrow 0.
$$
Thus $\Omega^1_{H_1}(\mathrm{log} E_{\mid H_1})\cong \mathcal O_{H_1}^{\oplus dim(H_1)}$ and $E_{H_1}\neq \emptyset$. By replacing $X$ with an appropriate open subset we may assume that $\T_{H_1}(-\mathrm{log} E_{H_1})\cong \mathcal O^{\oplus(n-1)}_{H_1}$. Let $G_2,\dots G_{n-2}$ be general hyperplanes and $H_i$, $T$, $S$ and $C$ as before, proceeding by induction we see that $\T_T(-\mathrm{log} C)\cong \mathcal O_{T}^{\oplus 2}$ and thus that $\T_S$ is locally free of rank $2$. \eproof\\

\textit{Proof of Theorem 1.1.}
Using the notation of Theorem 3.5 we assume that $Sing(X)\neq \emptyset$. With Lemma 3.2 (by induction) we can deduce that the surface $S$ constructed above necessarily has an isolated singularity. However $S$ is a log canonical surface with locally free tangent sheaf and thus smooth by Theorem 2.3. This contradicts the assumption.\eproof


\begin{thebibliography}{KMM87}

\bibitem[Che96]{Che96}
{\sc I.~A.Cheltsov}: \emph{Singularities of 3-Dimensional Varieties Admitting
an Ample Effective Divisor of Kodaira Dimension Zero}, Mathmatical Notes,
  Vol. 163, no. 4, 1970. {\sf\scriptsize UDC 512.774.42}
  
\bibitem[Dru13]{Druel13}
{\sc S. Druel}: \emph{The {Z}ariski-{L}ipman conjecture for log
canonical spaces},
2014 in Bulletin of the London Mathematical
Society Advance Access, June 9, 2014. {\sf\scriptsize 10.1112/blms/bdu040}

\bibitem[EV92]{EV92}
{\sc H.~Esnault and E.~Viehweg}: \emph{Lectures on vanishing theorems}, DMV
  Seminar, vol.~20, Birkh\"auser Verlag, Basel, 1992. {\sf\scriptsize MR1193913
  (94a:14017)}

\bibitem[Fle88]{Flenner88}
{\sc H.~Flenner}: \emph{Extendability of differential forms on nonisolated
  singularities}, Invent. Math. \textbf{94} (1988), no.~2.317--326
  {\sf\scriptsize MR958835 (89j:14001)}

\bibitem[GK14]{GK14}
{\sc P.~Graf and S.~J.~Kov{\'a}cs}: \emph{An optimal extension theorem for 1-forms and the Zariski-Lipman conjecture}, Doc. Math. \textbf{19}, 815--830, 2014.

\bibitem[GKK10]{GKK10}
{\sc D.~Greb, S.~Kebekus and S.~J.~Kov{\'a}cs}: \emph{Extension theorems for differential forms, and {B}ogomolov-{S}ommese
  vanishing on log canonical varieties}, Compositio Math. \textbf{1} 193--219,2010.
  
\bibitem[GKKP11]{GKKP11}
{\sc D.~Greb, S.~Kebekus, S.~J.~Kov{\'a}cs and T. Peternell} \emph{Differential forms on log canonical spaces}, Publ. Math. Inst. Hautes \'Etudes Sci. \textbf{114}, 87--169, 2011. 
 
\bibitem[KM98]{KM98}
{\sc J. Koll{\'a}r and S. Mori}: \emph{Birational geometry of algebraic varieties}, Cambridge Tracts in Mathematics, vol. 134, Cambridge University Press, Cambridge, 1998.
 
\bibitem[Kol07]{Kol07}
{\sc J.~Koll{\'a}r}: \emph{Lectures on resolution of singularities}, Annals of
  Mathematics Studies, vol. 166, Princeton University Press, Princeton, NJ,
  2007. {\sf\scriptsize MR2289519}

\bibitem[Lip65]{Lip65}
{\sc J.~Lipman}: \emph{ Free derivation modules on algebraic varieties},
Amer. J. Math. \textbf{87} (1965), no. 4.

\bibitem[Rei80]{Reid80}
{\sc M.~Reid}: \emph{Canonical Threefolds}, G{\'e}om{\'e}trie Alg{\'e}brique Angers, A. Beauville ed., Sijthoff
and Noordhoff, 1980.

\bibitem[Rei87]{Reid87}
{\sc M.~Reid}: \emph{Young person's guide to canonical singularities},
  Algebraic geometry, Bowdoin, 1985 (Brunswick, Maine, 1985), Proc. Sympos.
  Pure Math., vol.~46, Amer. Math. Soc., Providence, RI, 1987.
  {\sf\scriptsize MR927963 (89b:14016)}

\bibitem[Sai80]{Saito80}
{\sc K. Saito}: \emph{Theory of logarithmic differential forms and logarithmic vector fields},
 Journal of the Faculty of Science, the University of Tokyo. Sect. 1 A, Mathematics 27 (2),
 \href{http://repository.dl.itc.u-tokyo.ac.jp/dspace/bitstream/2261/6265/1/jfs270202.pdf}{JAIRO:jfs270202}.

\bibitem[vSS85]{vSS85}
{\sc D.~van Straten and J.~Steenbrink}: \emph{Extendability of holomorphic
  differential forms near isolated hypersurface singularities}, Abh. Math. Sem.
  Univ. Hamburg \textbf{55} (1985). {\sf\scriptsize MR831521
  (87j:32025)}

\bibitem[Wah85]{Wah85}
{\sc J.~M.~Wahl}: \emph{A characterization of quasi-homogeneous Gorenstein surface singularities}, Compositio Math. 55 (1985), no. 3, 269–288.


  
\end{thebibliography}
\end{document}